\documentclass{ws-procs9x6}
\usepackage{siunitx}
\usepackage{xargs}
\usepackage[capitalise]{cleveref}
\usepackage{subfigure}
\usepackage{mathtools}
\usepackage{booktabs}
\DeclareMathAlphabet{\mathpzc}{OT1}{pzc}{m}{it}

\usepackage{macros}

\usepackage[mode=buildnew,subpreambles=true]{standalone}
\usepackage{shellesc}
\usepackage{tikz,tikzscale}
\usetikzlibrary{arrows,calc,fit,patterns,positioning,spy,backgrounds,decorations.fractals,shapes.misc}
\tikzset{boximg/.style={remember picture,black,thick,draw,inner sep=0pt,outer sep=0pt}}
\usepackage{pgf}
\usepackage{pgfplots}
\usepgfplotslibrary{statistics}
\usepgfplotslibrary{fillbetween}
\usepgfplotslibrary{colormaps}
\pgfplotsset{compat=1.16}
\pgfplotsset{lua backend=true}
\usepgfplotslibrary{groupplots,dateplot}

\begin{document}

\title{Entropy stable shock capturing for high-order DGSEM on moving meshes}
\author{Anna Schwarz$^{a,*}$, Jens Keim$^{a}$, Christian Rohde$^{b}$, and Andrea Beck$^{a}$}
\address{$^a$Institute of Aerodynamics and Gas Dynamics, University of Stuttgart \\
         $^b$Institute of Applied Analysis and Numerical Simulation, University of Stuttgart}

\begin{abstract}
In this paper, a shock capturing for high-order entropy stable discontinuous Galerkin spectral element methods on moving
meshes is proposed using
Gauss--Lobatto nodes. The shock capturing is achieved via the convex blending of the high-order scheme with a low-order finite
volume subcell operator. The free-stream and convergence properties of the hybrid scheme are demonstrated numerically along with the entropy stability and shock capturing capabilities.
\end{abstract}

\bodymatter

\section{Introduction}
The numerical approximation of conservation laws on time-dependent domains with moving boundaries is an interesting but challenging field
of research with various applications such as shock tracking~\cite{Shi2022} or fluid-structure interactions~\cite{Donea1982}.
An example of a moving mesh approach is the arbitrary Lagrangian--Eulerian (ALE) method~\cite{Farhat2001} which theoretically allows for arbitrary mesh movements, limited only by the resulting mesh quality.
In the ALE ansatz, the time-dependent domain is transformed to a stationary, i.e., time-independent reference domain which introduces a mesh velocity in the convective flux of the conservation law.
Although this formulation is at first glance non-conservative, a conservative form can be retrieved by fulfilling an
additional conservation law denoted as the geometric conservation law (GCL).

High-order methods such as the discontinuous Galerkin (DG) spectral element method (DGSEM) are well suited for the accurate
approximation of complex, time-dependent multi-scale problems
due to their flexibility as well as low numerical dispersion and dissipation compared to low-order
schemes~\cite{Ainsworth2004a}.
DGSEM is a special class of DG methods in which the nodal solution is approximated by tensor-product Lagrange polynomials and the interpolation and quadrature nodes are collocated.
However, the application of high-order methods to nonlinear hyperbolic equations is inherently prone to (a) instabilities due
to aliasing errors and (b) Gibbs oscillations in the vicinity of strong discontinuities.
To alleviate these problems, additional stabilisation techniques are required, usually in combination with appropriate
shock capturing procedures.
Common stabilisation techniques include e.g. slope limiting~\cite{Krivodonova2007a}, filtering~\cite{Hesthaven2008}, polynomial de-aliasing through over-integration~\cite{Kirby2003} or entropy stable schemes~\cite{Fisher2013,Carpenter2014}.
The shock capturing may be achieved via a locally $h$-refined low-order finite volume (FV) scheme~\cite{Sonntag2017}, artificial viscosity~\cite{Persson2006}, flux correction \cite{Kuzmin2005}, or a convex blending of a high-order with a lower-order solution~\cite{Hennemann2021}, only to mention some possibilities.

This paper focuses on the construction of an entropy stable hybrid DGSE/FV method on moving meshes using the ALE approach and a
convex blending of a high-order DGSEM with a low-order FV scheme.
For the construction of the hybrid scheme, ideas from the entropy stable moving mesh DGSEM discussed by Schn{\"{u}}cke et al.~\cite{Schnucke2018}, and the entropy
stable hybrid DGSE/FV operator on static meshes developed by Hennemann et al.~\cite{Hennemann2021} are employed.

The outline is as follows: In~\cref{sec:ale}, the ALE formulation of a hyperbolic conservation system is briefly introduced. The hybrid DGSE/FV operator is presented in~\cref{sec:numerics}, followed in~\cref{sec:validation} by a validation of the resulting scheme. This paper closes by a short conclusion in~\cref{sec:conclusion}.

\section{Arbitrary Lagrangian-Eulerian formulation}
\label{sec:ale}
In this work, a compressible, inviscid fluid is considered which is governed by the Euler equations,
written in conservative form as
\begin{align}
  \cons_t(\pos,t) + \divX \physflux(\cons(\pos,t)) = \Null,
  \label{eq:theory:NSE}
\end{align}
and defined on the time-varying domain $\Omega_t \subset \Rd, \ t \in \R^+_0$ with dimension $d=3$.
\Cref{eq:theory:NSE} can be transformed to a time-independent reference domain, $E = [-1,1]^d$ with $\cons(\refpos,\tau)$, by
introducing a differentiable mapping $\boldsymbol{\chi}: E \to
\Omega_t, (\refpos,\tau) \mapsto (\pos,t)$, with time $\tau \in \R^+_0$.
In~\cref{eq:theory:NSE}, $\cons=\arr{\rho, \rho \vel, \rho e} \in \R^{\nq}$ with $\nq=5$, denotes the vector of the conserved variables, consisting
of the density $\rho$, the velocity vector $\vel=\arr{\vel[1],\vel[2],\vel[3]} \in \R^3$ , and the total energy $e$ per unit volume,
and $\fphysref \in \R^{3 \times \nq}$ are the contravariant fluxes, comprised of the well-known convective fluxes $\fphys$.
The ALE formulation of~\cref{eq:theory:NSE} is
\begin{align}
  \lr{\J \cons}_t + \gradientXI \cdot \fphysref  = \Null, \hspace{0.5cm} \fphysref = \M^\transpose \lr{\fphys - \cons \otimes \meshVel},
  \label{eq:theory:nse}
\end{align}
with the mesh velocity, $\meshVel = \lr{\boldsymbol{\chi}}_\tau \in \R^3$, and the geometric conservation law
\begin{align}
  \J_t &= \divXI \lre{\lr{\mathrm{adj}~\Jm^\transpose} \cdot \meshVel}.
  \label{eq:theory:gcl}
\end{align}
Here, $\tau = t$ (just a change of the notation) and the contravariant fluxes are $\fphysref$ including the ALE flux, $\cons \meshVel$.
In~\cref{eq:theory:nse}, the matrix $\M=\mathrm{adj} \Jm^\transpose \otimes \unit_{\nq}$ is defined as the adjoint of the
Jacobian matrix $\Jm$ of the mapping and $\J=\mathrm{det} \Jm$ is the determinant of the Jacobian.
Since interpolation and differentiation only commute if the discretization error is zero, the metric terms, $\M^\transpose$,
are approximated by the conservative curl form~\cite{Kopriva2006} to guarantee that the metric identities, $\divXI \Jm = 0$, hold on the discrete level.
The reader is referred to~\cite{Minoli2011,Schnucke2018} for a more detailed derivation of the ALE formulation.
The Euler equations are closed by the equation of state of a calorically perfect gas.

\section{Hybrid DG/FV operator on moving meshes}
\label{sec:numerics}
Following the method of lines approach,~\cref{eq:theory:nse,eq:theory:gcl} are discretized in space by the hybrid entropy stable DGSE/FV
operator on Gauss--Lobatto nodes and advanced in time by a high-order low-storage Runge--Kutta RK(4,5) method~\cite{Carpenter1994}.
The solution is approximated by the tensor-product of 1D Lagrange polynomials $\ell(\refpos)$ of order $\ppn$, resulting in 1D in
$\cons\approx\consh=\inte(\cons) : = \sum_{i=1}^\ppn \hat{\cons}(t) \ell_i(\refpos)$.
An entropy stable DGSE formulation relies on the fact that Fisher et al.~\cite{Fisher2013} have shown that a differentiation operator which fulfills the diagonal-norm summation-by-parts
(SBP) property such as the DGSEM on Gauss--Lobatto nodes can be reformulated in a flux differencing form that resembles a low-order finite volume (FV) subcell scheme.
According to Tadmor~\cite{Tadmor1987}, this enables the derivation of a discrete entropy condition, given an adequate two-point flux.
The entropy stable DGSE can be retrieved using the diagonal-norm SBP property, ${\Qm}^T+\Qm=\Bm$, and an adequate two-point flux
$\fsharp$, with $\Qm =  \Mm \Dm, \ \Bm =
\mathrm{diag}(-1,0,\ldots,0,1)$, the polynomial derivative matrix $\Dm_{ij}=\ell_j'(\refpos_i)$, and the mass matrix $\Mm =
\mathrm{diag}(\w_0,\ldots,\w_{\mathcal{N}})$ comprised of the quadrature weights $\w$.
Following the convex blending approach in Hennemann et al.~\cite{Hennemann2021} with the element-wise convex blending coefficient
$\FValpha \in [0,1]$, the hybrid DGSE/FV operator can be written as
\begin{equation*}
\resizebox{.96\hsize}{!}{\ensuremath{
  \Mm(\Jac \consh)_t = (1-\FValpha) \sum\limits_{k=1}^3 \left[- (2\Qm \circ \fcontsharp_k) \mathbf{1}
    - \Bm \begin{pmatrix} \fcont^\ast_{L,k} - \fcontsharp_{L,k} \\
  \fcont^\ast_{R,k} - \fcontsharp_{R,k} \end{pmatrix}\right]
  + \FValpha \sum\limits_{k=1}^3 \Bm^{FV} \lre{ \begin{matrix} \fcont^{FV}_{L,k} \\
  \fcont^{FV}_{R,k} \end{matrix}},
}}
\end{equation*}
and the discrete geometric conservation law as
\begin{equation*}
\resizebox{.96\hsize}{!}{\ensuremath{
        \Mm (\Jac)_t = (1-\FValpha) \sum\limits_{k=1}^3 \lre{(2\Qm \circ \meshvelsharp) \mathbf{1}
          + \Bm \begin{pmatrix} \fphysgclref^\ast_{L,k} - \meshvelsharp_{L,k} \\ \fphysgclref^\ast_{R,k} - \meshvelsharp_{R,k} \end{pmatrix}}
                       - \FValpha \Bm^{FV} \sum\limits_{k=1}^3 \lre{ \begin{matrix} {\fphysgclref}^{FV}_{L,k} \\
                      {\fphysgclref}^{FV}_{R,k} \end{matrix}},
}}
\end{equation*}
where $\Bm^{FV}=\mathrm{diag}([1,-1])$, see Hennemann et al.~\cite{Hennemann2021} for further details.
The contravariant fluxes are
\begin{equation}
  \fcontsharp_{ij} = \avg{\M^\transpose}_{ij} \lr{\fsharp_{ij} - \avg{\meshVel}_{ij} \qsharp_{ij}} \hspace{0.5cm} \meshvelsharp =\avg{\M^\transpose} \avg{\fphysgcl}
\end{equation}
with the averaging operator $\avg{\cdot}_{ij}=[(\cdot)_i+(\cdot)_j]/2$, an entropy conservative flux $\fsharp$ and an adequate state function
$\qsharp$~\cite{Schnucke2018}.
Neighboring cells are weakly coupled via the numerical flux $f^\ast$ normal to the face,
\begin{equation}
  \fcont^\ast_{L} = f^\ast(\cons_{L},\cons^+;\meshVel_{L}, \normalvec) \J_f, \hspace{0.5cm}
  \fphysgclref^\ast_{L} = (\meshVel_{L} \cdot \normalvec) \J_f,
\end{equation}
where $\J_f = \abs{\mathrm{adj}{\Jm_k}^\transpose \cdot {\normalvecref}_k}$, $\cons_{L}$ and $\cons^+$ are the values of $\consh$ right and left to the left element interface and $\normalvecref$ is the normal vector in reference space.
The numerical fluxes for the FV scheme are given as
\begin{equation}
  \fcont^{FV}_{L,k} = f^\ast(\cons^{FV}_{+,k},\cons^{FV}_{L,k}; \meshVel^{FV}_{L,k},\normalvec^{FV}_k)\J_f^{FV}, \hspace{0.5cm}
  {\fphysgclref}^{FV}_{L,k} = (\meshVel_{L,k}^{FV} \cdot \normalvec^{FV}_k)\J_f^{FV}
    ,
\end{equation}
with $\J_f^{FV} = \abs{\mathrm{adj}{\Jm^{FV}_k}^\transpose \cdot {\normalvecref}_k^{FV}}$, the mesh velocities $\meshVel^{FV}_k$ and the unit normal $\normalvecref^{FV}_k$ and physical normal vector $\normalvec^{FV}_k$ at the subcell
interface.
To ensure that the elements remain non-overlapping for $t>0$, the mesh velocities have to be continuous across the cells and thus
unique.
It has to be noted that compared to Hennemann et al.~\cite{Hennemann2021}, the surface fluxes have to be also blended since in
general $\meshVel_f^{FV} \neq \meshVel_f$. The same applies if a second-order FV subcell scheme is used, where
$\cons_f^{FV}\neq\cons_f$. To ensure conservation, a unique blending coefficient is chosen on the sides, $\alpha_f =
\max(\alpha_L,\alpha+)$.
Combining the proofs in Hennemann et al.~\cite{Hennemann2021} and Schn{\"{u}}cke et al.~\cite{Schnucke2018}, it can be shown that for an adequate two-point flux and a first-order FV method this hybrid scheme is either entropy conservative, under the condition that $\FValpha \in
[0,1]$, see~\cite{SchwarzPHD}. A numerical proof is given below.
The open-source framework FLEXI\footnote{www.flexi-project.org}~\cite{Krais2021,Kempf2024} is used as a solver which includes the numerical methods mentioned below.

\section{Numerical Experiments}
\label{sec:validation}

For the following surveys, if not stated otherwise, the computational domain is of size $\Omega_t = [-1,1]^3$ with periodic boundary conditions and
a convective time step restriction of $\mathrm{CFL}=0.9$.
The mesh was deformed according to Minoli et al.~\cite{Minoli2011}, where each corner node was displaced sinusoidally in time.
The numerical flux function $f^\ast$ is approximated by the Rusanov flux function and the entropy conservative flux function
of Chandrashekar et al.~\cite{Chandrashekar2013} (CH) is employed as two-point flux.

\subsection{Free-stream preservation}
For a static mesh, conservation is guaranteed if the discrete metric identities are satisfied, while for a
moving mesh, additionally the discrete geometric conservation law has to be taken into account, cf.~\cref{sec:ale}.
The free-stream preservation property, i.e., the ability to preserve a constant state over time, was evaluated for an arbitrary constant state, here given as $\rho=1$, $\vel=0.3$ and
$p=17.857$, which was advanced in time until $t=1$ with a randomly distributed $\FValpha\in[0,1)$.
The computational domain was discretized by $4^3$ elements.
The results in~\cref{tab:validation:freestream} highlight the free-stream preservation property of the hybrid DG/FV subcell scheme on moving meshes.
\begin{table}
  \centering
  \resizebox{\textwidth}{!}{\begin{tabular}{c|ccccc}
          & $\rho$ & $\rho \vel[1]$ & $\rho \vel[2]$ & $\rho \vel[3]$ & $\rho e$   \\ \hline
      $\ppn=3$  & 2.350883E-16 & 2.144119E-16 & 1.693604E-16 & 1.719768E-16 & 1.013106E-14 \\
      $\ppn=4$  & 2.133170E-16 & 1.051820E-16 & 1.026588E-16 & 9.937051E-17 & 1.072318E-14 \\
      $\ppn=5$  & 1.998885E-16 & 9.198981E-17 & 8.576905E-17 & 8.639587E-17 & 9.002652E-15 \\\hline
  \end{tabular}}
  \caption{Free-stream preservation: $L_2$ error of $\cons$.}
  {\vspace*{0.2cm}\footnotesize Table 1: Free-stream preservation: $L_2$-error of $\cons$.}
  \label{tab:validation:freestream}
\end{table}

\subsection{Convergence rates}

\begin{figure}[htbp]
    \begin{tikzpicture}
    \setlength{\fboxsep}{0pt}

    \small
    \begin{groupplot}[
          group style={
            group size=2 by 1,
            horizontal sep={0.17\linewidth},
            vertical sep={0.15\linewidth}
          },
          width=0.9\linewidth,
      ]

    \def\figwidth{0.48\linewidth}

    \nextgroupplot[
      enlargelimits=false,
      width=0.45\columnwidth,height=0.45\columnwidth,
      xlabel={$\mathcal{N}$},ylabel={$L_2$},
      xmin=0.6, xmax=10.4,
      ymin=1.e-9, ymax=1.e-1,
      ymode=log,
      axis y line*=left, axis x line*=bottom,
    ]
    \addplot [semithick, black, mark=o, mark size=2, mark options={solid}]
    table {%
    1  3.40e-02
    2  2.80e-03
    3  2.41e-04
    4  1.83e-05
    5  1.18e-06
    6  6.90e-08
    7  1.81e-08
    8  1.45e-08
    9  1.21e-08
    10 1.02e-08
    };
    \addlegendentry{DG}

    \addplot [semithick, black, mark=square, mark size=2, mark options={solid}]
    table {%
    1  3.40e-02
    2  2.80e-03
    3  2.41e-04
    4  1.83e-05
    5  1.18e-06
    6  6.89e-08
    7  1.80e-08
    8  1.44e-08
    9  1.20e-08
    10 1.01e-08
    };
    \addlegendentry{DG $\nu=0$}

    \nextgroupplot[
      enlargelimits=false,
      width=0.45\columnwidth,height=0.45\columnwidth,
      legend cell align={left},
      legend style={draw=none,fill opacity=0.0, draw opacity=1, text opacity=1, at={(0.7,0.97)}},
      xlabel={$\Delta x$},ylabel={$L_2$},
      xmin=0.0512656307826354, xmax=2.10956947206785,
      xmode=log,ymode=log,
      ymin=1.21074290991898e-10, ymax=5.E-1,
      axis y line*=left, axis x line*=bottom,
    ]
    \addplot [semithick, black, mark=square, mark size=2, mark options={solid}, forget plot]
    table {%
    0.125 2.10e-08
    0.25  5.82e-07
    0.5   1.83e-05
    1     5.48e-04
    2     1.34e-02
    };
    \addplot [semithick, black, mark=o, mark size=2, mark options={solid}, forget plot]
    table {%
    0.125  2.10e-08
    0.25   5.82e-07
    0.5    1.83e-05
    1      5.48e-04
    2      1.34e-02
    };
    \addplot [semithick, black, dashed]
    table {%
    0.0625 9.5367431640625e-09
    0.125  3.0517578125e-07
    0.25   9.765625e-06
    0.5    0.0003125
    1      0.01
    2      0.32
    };
    \addlegendentry{$\propto \Delta x^{\mathcal{N}+1}$}

    \end{groupplot}

    \end{tikzpicture}
  \caption{Validation of the spatial discretization with curved faces for a static (DG ($\meshVel=0$)) and sinusoidally deformed
    domain (DG). Left: $p$-convergence on a $4^3$ grid with $\ppn \in [1,10]$. Right: $h$-convergence for $\ppn=4$. Grid sequence
  ranges from $2^3$ up to $32^3$.}
  \label{fig:validation:conv}
\end{figure}
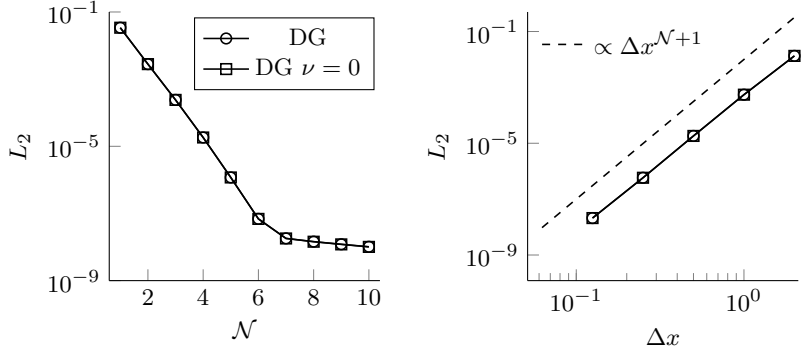

In the following, the $h$- and $p$-convergence properties of the DGSE operator (without FV) on a moving domain are investigated and compared to the results on a static domain.
In the $p$-convergence survey, the grid size was fixed to $4^3$ elements in each direction and the polynomial degree was varied, while
for the $h$-convergence study, $\Omega_t$ was discretized using $4^3$ to $32^3$ elements in each direction and $\ppn=4$.
The chosen error norm is the discrete $L_2$ error of the density.
The results in~\cref{fig:validation:conv} highlight the spatial convergence properties of the ALE
DGSEM for moving and static meshes using $\mathrm{CFL}=0.1$.

\subsection{Entropy conservation and stability}
\begin{figure}[htbp]
  \centering
  \includegraphics[]{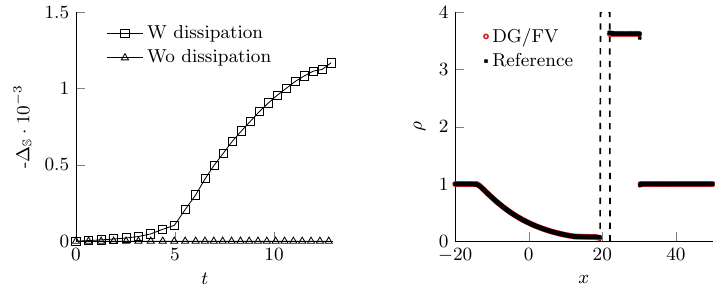}
  \caption{Left: Temporal evolution of the integral entropy conservation
    errors for the Euler equations using the TGV as a test case with and without additional surface dissipation. Right: Density
    distribution at $t=11$ for the moving piston test case. The dashed lines highlight the position of the piston. \label{fig:tgv_piston}}
\end{figure}

For the investigation of the entropy or the kinetic energy conservation property, a common numerical test case is the inviscid
Taylor-Green vortex~\cite{Shu2005a} at $M=0.1$.
The computational domain was of size $\Omega_t = [0,2\pi]^3$ with $16^3$ elements, $\ppn=3$ and periodic boundary conditions.
A time step restriction of $CFL=0.25$, a constant blending coefficient of $\FValpha=0.3$ and a first-order FV
subcell operator were chosen.
To further illustrate the ability of the framework to deal with curvilinear element representations, the mesh was
deformed according to a standing wave in 3D~\cite{Schnucke2018} and represented by a tensor-product of 1D Lagrange polynomials of order $\ppngeo=2$.
The temporal evolution of the integral entropy error $\Delta_{\mathbb{S}}$, defined in~\cite{Schnucke2018}, on the left
of~\cref{fig:tgv_piston} illustrates that the CH
flux recovers the entropy conservation property for the case without additional surface dissipation (only the central part of the
Riemann solver was employed), otherwise the scheme is entropy stable.

\subsection{Shock capturing}
To test the shock capturing, the 1D piston shock tube problem~\cite{Tan2011} is considered, consisting of a right moving shock wave and a left travelling rarefaction wave.
The piston moves at $M=2$ in an initially quiescent fluid with $\arr{\rho, \vel, p} = \arr{1,\mathbf{0},5/7}$.
The computational domain $\Omega = [-20,50] \times [0,10]^2$ is discretized using $280 \times 1^2$ elements in each direction with periodic boundary conditions in $y$- and $z$- direction,
while Dirichlet boundary conditions with the initial state are imposed in $x$-direction. A second-order FV subcell operator was
chosen.
The predicted solution is in good agreement with the reference solution~\cite{Liepmann1959}, cf.~\cref{fig:tgv_piston} (right).

\section{Conclusion}
\label{sec:conclusion}
In this work, a hybrid entropy stable DGSE/FV operator for moving meshes using Gauss--Lobatto nodes was proposed.
It was shown that the scheme is free-stream preserving and retains the convergence properties of the hybrid scheme on static
meshes. Finally, the entropy stability and shock capturing capabilities were demonstrated.

\section*{Acknowledgements}
This work was funded by the European Union and has received funding from the European High Performance Computing Joint Undertaking (JU)
and Sweden, Germany, Spain, Greece, and Denmark under grant agreement No 101093393.
Moreover, this research was funded by the DFG under Germany's Excellence Strategy EXC 2075-390740016.
Further we want to gratefully acknowledge funding by the DFG through SPP 2410 Hyperbolic Balance Laws in Fluid Mechanics: Complexity, Scales, Randomness (CoScaRa).

\bibliographystyle{ws-procs9x6}
\bibliography{references.bib}

\begin{thebibliography}{10}

\bibitem{Shi2022}
A.~Shi, P.-O. Persson and M.~Zahr, {Implicit shock tracking for unsteady flows
  by the method of lines}, {\em J. of Comput. Phys.} {\bf 454}, p. 110906
  (2022).

\bibitem{Donea1982}
J.~Donea, S.~Giuliani and J.~Halleux, {An arbitrary Lagrangian-Eulerian finite
  element method for transient dynamic fluid-structure interactions}, {\em
  Comput. Methods Appl. Mech. Eng.} {\bf 33}, 689  (1982).

\bibitem{Farhat2001}
C.~Farhat, P.~Geuzaine and C.~Grandmont, {The discrete geometric conservation
  law and the nonlinear stability of ALE schemes for the solution of flow
  problems on moving grids}, {\em J. of Comput. Phys.} {\bf 174}, 669  (2001).

\bibitem{Ainsworth2004a}
M.~Ainsworth, {Dispersive and dissipative behaviour of high order discontinuous
  Galerkin finite element methods}, {\em J. of Comput. Phys.} {\bf 198}, 106
  (2004).

\bibitem{Krivodonova2007a}
L.~Krivodonova, {Limiters for high-order discontinuous Galerkin methods}, {\em
  J. Comput. Phys.} {\bf 226}, 879  (2007).

\bibitem{Hesthaven2008}
J.~S. Hesthaven and R.~M. Kirby, {Filtering in Legendre spectral methods}, {\em
  Mathematics of Computation} {\bf 77}, 1425  (2008).

\bibitem{Kirby2003}
R.~Kirby and G.~Karniadakis, {De-aliasing on non-uniform grids: Algorithms and
  applications}, {\em J. Comput. Phys.} {\bf 191}, 249 (nov 2003).

\bibitem{Fisher2013}
T.~C. Fisher and M.~H. Carpenter, High-order entropy stable finite difference
  schemes for nonlinear conservation laws: Finite domains, {\em J. Comput.
  Phys.} {\bf 252}, 518  (2013).

\bibitem{Carpenter2014}
M.~H. Carpenter, T.~C. Fisher, E.~J. Nielsen and S.~H. Frankel, {Entropy Stable
  Spectral Collocation Schemes for the Navier--Stokes Equations: Discontinuous
  Interfaces}, {\em SIAM J. Sci. Comput.} {\bf 36}, B835  (2014).

\bibitem{Sonntag2017}
M.~Sonntag and C.-D. Munz, {Efficient Parallelization of a Shock Capturing for
  Discontinuous Galerkin Methods using Finite Volume Sub-cells}, {\em J. Sci.
  Comput.} {\bf 70}, 1262  (2017).

\bibitem{Persson2006}
P.-O. Persson and J.~Peraire, Sub-cell shock capturing for discontinuous
  {G}alerkin methods, in {\em 44th AIAA Aerospace Sciences Meeting and
  Exhibit\/}, 2006.

\bibitem{Kuzmin2005}
D.~Kuzmin and M.~M{\"o}ller, {\em Algebraic Flux Correction II. Compressible
  Euler Equations}, in {\em Flux-Corrected Transport: Principles, Algorithms,
  and Applications\/},  (Springer, 2005), pp. 207--250.

\bibitem{Hennemann2021}
S.~Hennemann, A.~M. Rueda-Ram{\'{i}}rez, F.~J. Hindenlang and G.~J. Gassner, {A
  provably entropy stable subcell shock capturing approach for high order split
  form DG for the compressible Euler equations}, {\em J. Comput. Phys.} {\bf
  426}, p. 109935  (2021).

\bibitem{Schnucke2018}
G.~Schn{\"{u}}cke, N.~Krais, T.~Bolemann and G.~J. Gassner, {Entropy Stable
  Discontinuous Galerkin Schemes on Moving Meshes for Hyperbolic Conservation
  Laws}, {\em J. Sci. Comput.} {\bf 82}, p.~69  (2020).

\bibitem{Kopriva2006}
D.~A. Kopriva, Metric identities and the discontinuous spectral element method
  on curvilinear meshes, {\em J. of Sci. Comput.} {\bf 26}, p. 301  (2006).

\bibitem{Minoli2011}
C.~A.~A. Minoli and D.~A. Kopriva, Discontinuous {Galerkin} spectral element
  approximations on moving meshes, {\em J. Comput. Phys.} {\bf 230}, 1876
  (2011).

\bibitem{Carpenter1994}
M.~H. Carpenter and A.~Kennedy, {Fourth-Order 2N-Storage Runge-Kutta Schemes},
  {\em Nasa Technical Memorandum} {\bf 109112}, 1  (1994).

\bibitem{Tadmor1987}
E.~Tadmor, The numerical viscosity of entropy stable schemes for systems of
  conservation laws. {I}, {\em Math. Comput.} {\bf 49}, 91  (1987).

\bibitem{SchwarzPHD}
A.~Schwarz, {High-fidelity particle tracking and impact-induced deformations},
  PhD thesis, Universit{\"{a}}t Stuttgart2024.

\bibitem{Krais2021}
N.~Krais, A.~Beck, T.~Bolemann, H.~Frank, D.~Flad, G.~Gassner, F.~Hindenlang,
  M.~Hoffmann, T.~Kuhn, M.~Sonntag and C.-D. Munz, {FLEXI}: A high order
  discontinuous {G}alerkin framework for hyperbolic{\textendash}parabolic
  conservation laws, {\em Comput. Math. with Appl.} {\bf 81}, 186  (2021).

\bibitem{Kempf2024}
M.~Kurz, D.~Kempf, M.~P. Blind, P.~Kopper, P.~Offenh{\"{a}}user, A.~Schwarz,
  S.~Starr, J.~Keim and A.~Beck, {GAL{\AE}XI: Solving complex compressible
  flows with high-order discontinuous Galerkin methods on accelerator-based
  systems}, {\em Comput. Phys. Commun.} {\bf 306}, p. 109388  (2025).

\bibitem{Chandrashekar2013}
P.~Chandrashekar, Kinetic energy preserving and entropy stable finite volume
  schemes for compressible {Euler} and {Navier}-{Stokes} equations, {\em
  Commun. Comput. Phys.} {\bf 14}, 1252  (2013).

\bibitem{Shu2005a}
C.-W. Shu, W.-S. Don, D.~Gottlieb, O.~Schilling and L.~Jameson, {Numerical
  Convergence Study of Nearly Incompressible, Inviscid Taylor–Green Vortex
  Flow}, {\em J. Sci. Comput.} {\bf 24}, 1  (2005).

\bibitem{Tan2011}
S.~Tan and C.~W. Shu, {A high order moving boundary treatment for compressible
  inviscid flows}, {\em J. Comput. Phys.} {\bf 230}, 6023  (2011).

\bibitem{Liepmann1959}
H.~W. Liepmann and A.~Roshko, {\em {Elements of gasdynamics}} (John Wiley {\&}
  Sons, 1959).

\end{thebibliography}

\end{document}